\documentclass{amsart}
\usepackage{amsfonts}

\newtheorem{theorem}{Theorem}
\theoremstyle{plain}

\newtheorem{corollary}{Corollary}

\newtheorem{lemma}{Lemma}

\newtheorem{proposition}{Proposition}
\newtheorem{remark}{Remark}

\numberwithin{equation}{section}
\input{tcilatex}

\begin{document}
\title[Inequalities in Normed Linear Spaces]{On Some Discrete Inequalities
in Normed Linear Spaces}
\author{Sever S. Dragomir}
\address{School of Computer Science and Mathematics\\
Victoria University\\
PO Box 14428, Melbourne City\\
Victoria 8001, Australia.}
\email{sever.dragomir@vu.edu.au}
\urladdr{http://rgmia.vu.edu.au/dragomir}
\date{4 August, 2005}
\subjclass[2000]{Primary 46B05, 26D15.}
\keywords{Inequalities in normed spaces, Semi-inner products, Analytic
inequalities.}

\begin{abstract}
Some sharp discrete inequalities in normed linear spaces are obtained. New
reverses of the generalised triangle inequality are also given.
\end{abstract}

\maketitle

\section{Introduction}

Let $\left( X,\left\Vert \cdot \right\Vert \right) $ be a normed linear
space over the real or complex number field $\mathbb{K}$. The mapping $%
f:X\rightarrow \mathbb{R}$, $f\left( x\right) =\frac{1}{2}\left\Vert
x\right\Vert ^{2}$ is obviously convex on $\mathbb{R}$ and then there exists
the following limits:%
\begin{align*}
\left\langle x,y\right\rangle _{i}& :=\lim_{t\rightarrow 0-}\frac{\left\Vert
y+tx\right\Vert ^{2}-\left\Vert y\right\Vert ^{2}}{2t}, \\
\left\langle x,y\right\rangle _{s}& :=\lim_{\tau \rightarrow 0-}\frac{%
\left\Vert y+\tau x\right\Vert ^{2}-\left\Vert y\right\Vert ^{2}}{2\tau }
\end{align*}%
for any two vectors in $X.$ The mapping $\left\langle \cdot ,\cdot
\right\rangle _{s}$ $\left( \left\langle \cdot ,\cdot \right\rangle
_{i}\right) $ will be called the \textit{superior semi-inner product }(%
\textit{inferior semi-inner product}) associated to the norm $\left\Vert
\cdot \right\Vert .$

The following fundamental calculus rules are valid for these semi-inner
products (see for instance \cite[p. 27--32]{D1}):%
\begin{align}
\left\langle x,x\right\rangle _{p}& =\left\Vert x\right\Vert ^{2}\quad \text{%
for }x\in X;  \label{1.1} \\
\left\langle \lambda x,y\right\rangle _{p}& =\lambda \left\langle
x,y\right\rangle _{p}\quad \text{for }\lambda \geq 0\text{ and }x,y\in X;
\label{1.2} \\
\left\langle x,\lambda y\right\rangle _{p}& =\lambda \left\langle
x,y\right\rangle _{p}\quad \text{for }\lambda \geq 0\text{ and }x,y\in X;
\label{1.3} \\
\left\langle \lambda x,y\right\rangle _{p}& =\lambda \left\langle
x,y\right\rangle _{q}\quad \text{for }\lambda \leq 0\text{ and }x,y\in X;
\label{1.4} \\
\left\langle \alpha x,\beta y\right\rangle _{p}& =\alpha \beta \left\langle
x,y\right\rangle _{p}\quad \text{for }\alpha ,\beta \in \mathbb{R}\text{
with }\alpha \beta \geq 0\text{ and }x,y\in X;  \label{1.5} \\
\left\langle -x,y\right\rangle _{p}& =\left\langle x,-y\right\rangle
_{p}=-\left\langle x,y\right\rangle _{q}\quad \text{for }x,y\in X;
\label{1.6}
\end{align}%
where $p,q\in \left\{ s,i\right\} $ and $p\neq q.$

The following inequality is valid:%
\begin{align}
\frac{\left\Vert y+tx\right\Vert ^{2}-\left\Vert y\right\Vert ^{2}}{2t}&
\geq \left\langle x,y\right\rangle _{s}\geq \left\langle x,y\right\rangle
_{i}  \label{1.7} \\
& \geq \frac{\left\Vert y+sx\right\Vert ^{2}-\left\Vert y\right\Vert ^{2}}{2s%
},  \notag
\end{align}%
for any $x,y\in X$ and $s<0<t.$

An important result is the following \textit{Schwarz inequality}:%
\begin{equation}
\left\vert \left\langle x,y\right\rangle _{p}\right\vert \leq \left\Vert
x\right\Vert \left\Vert y\right\Vert \quad \text{for each }x,y\in X.
\label{1.8}
\end{equation}%
Also, the following properties of sub(super)-additivity should be noted:%
\begin{equation}
\left\langle x_{1}+x_{2},y\right\rangle _{s\left( i\right) }\leq \left( \geq
\right) \left\langle x_{1},y\right\rangle _{s\left( i\right) }+\left\langle
x_{2},y\right\rangle _{s\left( i\right) }  \label{1.9}
\end{equation}%
for each $x_{1},x_{2},y\in X.$

Another important property of \textquotedblleft
quasi-linearity\textquotedblright\ holds as well:%
\begin{equation}
\left\langle \alpha x+y,x\right\rangle _{p}=\alpha \left\Vert x\right\Vert
^{2}+\left\langle y,x\right\rangle _{p}  \label{1.10}
\end{equation}%
for any $x,y\in X$ and $\alpha $ a real number, where $p=s$ or $p=i.$

Finally, we mention the continuity property:%
\begin{equation}
\left\vert \left\langle y+z,x\right\rangle _{p}-\left\langle
z,x\right\rangle _{p}\right\vert \leq \left\Vert y\right\Vert \left\Vert
x\right\Vert  \label{1.11}
\end{equation}%
for each $x,y,z\in X$ and $p=s$ or $p=i.$

One of the most used inequalities in normed spaces is the triangle
inequality for several vectors, i.e.,%
\begin{equation}
\left\Vert \sum_{j=1}^{n}x_{j}\right\Vert \leq \sum_{j=1}^{n}\left\Vert
x_{j}\right\Vert  \label{1.12}
\end{equation}%
for any $x_{j}\in X,$ $j\in \left\{ 1,\dots ,n\right\} .$

The main aim of this paper is to point out some inequalities for norms of
the vectors $x_{j}$ and $\sum_{j=1}^{n}x_{j},$ including some reverses of
the triangle inequality in the multiplicative form, i.e., lower bounds for
the quantity%
\begin{equation*}
\frac{\left\Vert \sum_{j=1}^{n}x_{j}\right\Vert }{\sum_{j=1}^{n}\left\Vert
x_{j}\right\Vert },
\end{equation*}%
provided that not all $x_{j}$ are zero and satisfy some appropriate
conditions.

For classical results related to the reverse of the triangle inequality in
normed spaces see \cite{DM}, \cite{K}, \cite{DPF} and \cite{M}. For more
recent results, see \cite{D2}, \cite{D3}, \cite{AM1} and \cite{AM2}.

\section{The Results}

The following lemma is of interest itself as well.

\begin{lemma}
\label{l2.1}Let $\left( X,\left\Vert \cdot \right\Vert \right) $ be a normed
linear space. If $x,a\in X,$ then%
\begin{equation}
\left\langle x,a\right\rangle _{i}\geq \frac{1}{2}\left( \left\Vert
a\right\Vert ^{2}-\left\Vert x-a\right\Vert ^{2}\right) .  \label{2.1}
\end{equation}%
If $\left\Vert a\right\Vert >\left\Vert x-a\right\Vert ,$ then the constant $%
\frac{1}{2}$ cannot be replaced by a larger quantity.
\end{lemma}

\begin{proof}
Utilising the semi-inner product properties, we have by $\left( \ref{1.7}%
\right) $ that%
\begin{equation*}
\left\langle x,a\right\rangle _{i}=\lim_{s\rightarrow 0-}\frac{\left\Vert
a+sx\right\Vert ^{2}-\left\Vert a\right\Vert ^{2}}{2s}\geq \frac{\left\Vert
a+\left( -1\right) x\right\Vert ^{2}-\left\Vert a\right\Vert ^{2}}{2\left(
-1\right) }=\frac{\left\Vert a\right\Vert ^{2}-\left\Vert x-a\right\Vert ^{2}%
}{2}
\end{equation*}%
and the inequality (\ref{2.1}) is proved.

Now, assume that $\left\Vert a\right\Vert >\left\Vert x-a\right\Vert $ and
there exists a $C>0$ with the property that%
\begin{equation}
\left\langle x,a\right\rangle _{i}\geq C\left( \left\Vert a\right\Vert
^{2}-\left\Vert x-a\right\Vert ^{2}\right) .  \label{2.2}
\end{equation}%
Obviously $a\neq 0,$ and if we choose $x=\varepsilon a,$ $\varepsilon \in
\left( 0,1\right) ,$ then $\left\Vert a\right\Vert >\left\Vert
x-a\right\Vert $ since $\left\Vert x-a\right\Vert =\left( 1-\varepsilon
\right) \left\Vert a\right\Vert .$ Replacing $x$ in (\ref{2.2}) we get%
\begin{equation*}
\varepsilon \left\Vert a\right\Vert ^{2}\geq C\left( \left\Vert a\right\Vert
^{2}-\left( 1-\varepsilon \right) ^{2}\left\Vert a\right\Vert ^{2}\right) 
\end{equation*}%
giving%
\begin{equation*}
\varepsilon \geq C\left( 2\varepsilon -\varepsilon ^{2}\right) ,
\end{equation*}%
for any $\varepsilon \in \left( 0,1\right) .$ This is in fact $1\geq C\left(
2-\varepsilon \right) $ and if we let $\varepsilon \rightarrow 0+,$ we get $%
C\geq \frac{1}{2}.$
\end{proof}

\begin{remark}
As a coarser, but maybe more useful inequality, we can state that%
\begin{equation}
\left\langle x,a\right\rangle _{i}\geq \frac{1}{2}\left\Vert x\right\Vert
\left( \left\Vert a\right\Vert -\left\Vert x-a\right\Vert \right) ,
\label{2.3}
\end{equation}%
provided $\left\Vert a\right\Vert \geq \left\Vert x-a\right\Vert .$

We observe that (\ref{2.3}) follows from (\ref{2.1}) since, for $\left\Vert
a\right\Vert \geq \left\Vert x-a\right\Vert ,$ the triangle inequality gives:%
\begin{align*}
\frac{1}{2}\left( \left\Vert a\right\Vert ^{2}-\left\Vert x-a\right\Vert
^{2}\right) & =\frac{1}{2}\left( \left\Vert a\right\Vert -\left\Vert
x-a\right\Vert \right) \left( \left\Vert a\right\Vert +\left\Vert
x-a\right\Vert \right)  \\
& \geq \frac{1}{2}\left( \left\Vert a\right\Vert -\left\Vert x-a\right\Vert
\right) \left\Vert x\right\Vert .
\end{align*}%
It is an open question whether the constant $\frac{1}{2}$ in (\ref{2.3}) is
sharp.
\end{remark}

The following result may be stated.

\begin{theorem}
\label{t2.1}Let $\left( X,\left\Vert \cdot \right\Vert \right) $ be a normed
space and $x_{j}\in X,$ $j\in \left\{ 1,\dots ,n\right\} $, $a\in
X\backslash \left\{ 0\right\} .$ Then for any $p_{j}\geq 0,j\in \left\{
1,\dots ,n\right\} $ with $\sum_{j=1}^{n}p_{j}=1$ we have%
\begin{equation}
\left\Vert \sum_{j=1}^{n}p_{j}x_{j}\right\Vert \left\Vert a\right\Vert +%
\frac{1}{2}\sum_{j=1}^{n}p_{j}\left\Vert x_{j}-a\right\Vert ^{2}\geq \frac{1%
}{2}\left\Vert a\right\Vert ^{2}.  \label{2.4}
\end{equation}%
The constant $\frac{1}{2}$ in the right hand side of (\ref{2.4}) is best
possible in the sense that it cannot be replaced by a larger quantity.
\end{theorem}

\begin{proof}
We apply Lemma \ref{l2.1} on stating that%
\begin{equation*}
\left\langle x_{j},a\right\rangle _{i}+\frac{1}{2}\left\Vert
x_{j}-a\right\Vert ^{2}\geq \frac{1}{2}\left\Vert a\right\Vert ^{2}
\end{equation*}%
for each $j\in \left\{ 1,\dots ,n\right\} .$

Multiplying with $p_{j}\geq 0$ and summing over $j$ from 1 to $n,$ we get%
\begin{equation}
\sum_{j=1}^{n}p_{j}\left\langle x_{j},a\right\rangle _{i}+\frac{1}{2}%
\sum_{j=1}^{n}p_{j}\left\Vert x_{j}-a\right\Vert ^{2}\geq \frac{1}{2}%
\left\Vert a\right\Vert ^{2}\sum_{j=1}^{n}p_{j}.  \label{2.5}
\end{equation}%
Utilising the superadditivity property of the semi-inner product $%
\left\langle \cdot ,\cdot \right\rangle _{i}$ in the first variable (see 
\cite[p. 29]{D1}) we have%
\begin{equation}
\left\langle \sum_{j=1}^{n}p_{j}x_{j},a\right\rangle _{i}\geq
\sum_{j=1}^{n}p_{j}\left\langle x_{j},a\right\rangle _{i}.  \label{2.6}
\end{equation}%
By the Schwarz inequality applied for $\sum_{j=1}^{n}p_{j}x_{j}$ and $a,$ we
also have%
\begin{equation}
\left\Vert \sum_{j=1}^{n}p_{j}x_{j}\right\Vert \left\Vert a\right\Vert \geq
\left\langle \sum_{j=1}^{n}p_{j}x_{j},a\right\rangle _{i}.  \label{2.7}
\end{equation}%
Therefore, by (\ref{2.5}) \ -- (\ref{2.7}) we deduce the desired inequality (%
\ref{2.4}).

Now assume that there exists a $D>0$ with the property that%
\begin{equation}
\left\Vert \sum_{j=1}^{n}p_{j}x_{j}\right\Vert \left\Vert a\right\Vert +%
\frac{1}{2}\sum_{j=1}^{n}p_{j}\left\Vert x_{j}-a\right\Vert ^{2}\geq
D\left\Vert a\right\Vert ^{2},  \label{2.8}
\end{equation}%
for any $n\geq 1,$ $x_{j}\in X,$ $p_{j}\geq 0$, $j\in \left\{ 1,\dots
,n\right\} $ with $\sum_{j=1}^{n}p_{j}=1$ and $a\in X\backslash \left\{
0\right\} .$

If in (\ref{2.8}) we choose $n=1,$ $p_{1}=1,$ $x_{1}=\varepsilon a,$ $%
\varepsilon \in \left( 0,1\right) ,$ then we get%
\begin{equation*}
\varepsilon \left\Vert a\right\Vert ^{2}+\frac{1}{2}\left( 1-\varepsilon
\right) ^{2}\left\Vert a\right\Vert ^{2}\geq D\left\Vert a\right\Vert ^{2},
\end{equation*}%
giving%
\begin{equation*}
\varepsilon +\frac{1}{2}\left( 1-\varepsilon \right) ^{2}\geq D,
\end{equation*}%
for any $\varepsilon \in \left( 0,1\right) .$ Letting $\varepsilon
\rightarrow 0+,$ we deduce $D\leq \frac{1}{2}$ and the proof is complete.
\end{proof}

The following result may be stated as well:

\begin{proposition}
\label{p2.1}Let $x_{j},a\in X$ with $a\neq 0$ and $\left\Vert
x_{j}-a\right\Vert \leq \left\Vert a\right\Vert $ for each $j\in \left\{
1,\dots ,n\right\} .$ Then for any $p_{j}\geq 0,j\in \left\{ 1,\dots
,n\right\} $ with $\sum_{j=1}^{n}p_{j}=1$ we have%
\begin{equation}
\left\Vert \sum_{j=1}^{n}p_{j}x_{j}\right\Vert \left\Vert a\right\Vert +%
\frac{1}{2}\sum_{j=1}^{n}p_{j}\left\Vert x_{j}\right\Vert \left\Vert
x_{j}-a\right\Vert \geq \frac{1}{2}\left\Vert a\right\Vert
\sum_{j=1}^{n}p_{j}\left\Vert x_{j}\right\Vert .  \label{2.9}
\end{equation}
\end{proposition}

\begin{proof}
From (\ref{2.3}) we have%
\begin{equation*}
\left\langle x_{j},a\right\rangle _{i}+\frac{1}{2}\left\Vert
x_{j}\right\Vert \left\Vert x_{j}-a\right\Vert \geq \frac{1}{2}\left\Vert
a\right\Vert \left\Vert x_{j}\right\Vert
\end{equation*}%
for any $j\in \left\{ 1,\dots ,n\right\} .$

The proof follows in the same manner as in Theorem \ref{t2.1} and we omit
the details.
\end{proof}

The following reverse of the generalised triangle inequality may be stated:

\begin{theorem}
\label{t2.2}Let $x_{j}\in X\backslash \left\{ 0\right\} $ and $a\in
X\backslash \left\{ 0\right\} $ such that $\left\Vert a\right\Vert \geq
\left\Vert x_{j}-a\right\Vert $ for each $j\in \left\{ 1,\dots ,n\right\} .$
Then for any $p_{j}\geq 0,j\in \left\{ 1,\dots ,n\right\} $ with $%
\sum_{j=1}^{n}p_{j}=1$ we have%
\begin{equation}
\frac{\left\Vert \sum_{j=1}^{n}p_{j}x_{j}\right\Vert }{\sum_{j=1}^{n}p_{j}%
\left\Vert x_{j}\right\Vert }\geq \frac{1}{2}\min_{1\leq j\leq n}\left\{ 
\frac{\left\Vert a\right\Vert ^{2}-\left\Vert a-x_{j}\right\Vert ^{2}}{%
\left\Vert x_{j}\right\Vert \left\Vert a\right\Vert }\right\} \left( \geq
0\right) .  \label{2.10}
\end{equation}%
The constant $\frac{1}{2}$ is best possible in (\ref{2.10}).
\end{theorem}

\begin{proof}
Let us denote%
\begin{equation*}
\rho :=\min_{1\leq j\leq n}\left\{ \frac{\left\Vert a\right\Vert
^{2}-\left\Vert a-x_{j}\right\Vert ^{2}}{\left\Vert x_{j}\right\Vert }%
\right\} .
\end{equation*}%
From Lemma \ref{l2.1} we have%
\begin{equation*}
\frac{\left\langle x_{j},a\right\rangle _{i}}{\left\Vert x_{j}\right\Vert }%
\geq \frac{1}{2}\cdot \frac{\left\Vert a\right\Vert ^{2}-\left\Vert
x_{j}-a\right\Vert ^{2}}{\left\Vert x_{j}\right\Vert }\geq \frac{1}{2}\rho 
\end{equation*}%
for each $j\in \left\{ 1,\dots ,n\right\} .$ Therefore%
\begin{equation*}
\left\langle x_{j},a\right\rangle _{i}\geq \frac{1}{2}\rho \left\Vert
x_{j}\right\Vert ,\quad j\in \left\{ 1,\dots ,n\right\} .
\end{equation*}%
Multiplying with $p_{j}$ and summing over $j$ from 1 to $n$ we obtain%
\begin{equation}
\sum_{j=1}^{n}p_{j}\left\langle x_{j},a\right\rangle _{i}\geq \frac{1}{2}%
\rho \sum_{j=1}^{n}p_{j}\left\Vert x_{j}\right\Vert ,  \label{2.11}
\end{equation}%
and since:%
\begin{equation}
\left\Vert \sum_{j=1}^{n}p_{j}x_{j}\right\Vert \left\Vert a\right\Vert \geq
\left\langle \sum_{j=1}^{n}p_{j}x_{j},a\right\rangle _{i}\geq
\sum_{j=1}^{n}p_{j}\left\langle x_{j},a\right\rangle _{i},  \label{2.12}
\end{equation}%
hence by (\ref{2.11}) and (\ref{2.12}) we deduce the desired result (\ref%
{2.10}).

Now, assume that there exists a constant $E>0$ such that%
\begin{equation}
\frac{\left\Vert \sum_{j=1}^{n}p_{j}x_{j}\right\Vert }{\sum_{j=1}^{n}p_{j}%
\left\Vert x_{j}\right\Vert }\geq E\cdot \min_{1\leq j\leq n}\left\{ \frac{%
\left\Vert a\right\Vert ^{2}-\left\Vert a-x_{j}\right\Vert ^{2}}{\left\Vert
x_{j}\right\Vert \left\Vert a\right\Vert }\right\} ,  \label{2.13}
\end{equation}%
provided $\left\Vert a\right\Vert \geq \left\Vert x_{j}-a\right\Vert ,$ $%
j\in \left\{ 1,\dots ,n\right\} .$

If we choose $x_{1}=\cdots =x_{n}=\varepsilon a,$ $\varepsilon \in \left(
0,1\right) ,$ and $p_{1}=...=p_{n}=\frac{1}{n}$, then we get%
\begin{equation*}
1\geq E\cdot \frac{\left\Vert a\right\Vert ^{2}-\left( 1-\varepsilon \right)
^{2}\left\Vert a\right\Vert ^{2}}{\varepsilon \left\Vert a\right\Vert ^{2}},
\end{equation*}%
giving%
\begin{equation*}
1\geq E\left( 2-\varepsilon \right) 
\end{equation*}%
for any $\varepsilon \in \left( 0,1\right) .$ Letting $\varepsilon
\rightarrow 0+,$ we deduce $E\leq \frac{1}{2}$ and the proof is complete.
\end{proof}

The following result may be stated as well:

\begin{proposition}
\label{p2.2}Let $x_{j},a\in X\backslash \left\{ 0\right\} ,$ $j\in \left\{
1,\dots ,n\right\} $ such that $\left\Vert x_{j}-a\right\Vert \leq
\left\Vert a\right\Vert .$ Then for any $p_{j}\geq 0,j\in \left\{ 1,\dots
,n\right\} $ with $\sum_{j=1}^{n}p_{j}=1$ we have%
\begin{equation}
\frac{\left\Vert \sum_{j=1}^{n}p_{j}x_{j}\right\Vert }{\sum_{j=1}^{n}p_{j}%
\left\Vert x_{j}\right\Vert }\geq \frac{\left( \left\Vert a\right\Vert
-\max_{1\leq j\leq n}\left\Vert x_{j}-a\right\Vert \right) }{2\left\Vert
a\right\Vert }\quad \left( \geq 0\right) .  \label{2.14}
\end{equation}
\end{proposition}

\begin{proof}
From (\ref{2.3}) we have%
\begin{align*}
\frac{\left\langle x_{j},a\right\rangle _{i}}{\left\Vert x_{j}\right\Vert }&
\geq \frac{1}{2}\left( \left\Vert a\right\Vert -\left\Vert
x_{j}-a\right\Vert \right) \\
& \geq \frac{1}{2}\min_{1\leq j\leq n}\left( \left\Vert a\right\Vert
-\left\Vert x_{j}-a\right\Vert \right) \\
& =\frac{1}{2}\left( \left\Vert a\right\Vert -\max_{1\leq j\leq n}\left\Vert
x_{j}-a\right\Vert \right) .
\end{align*}%
Now the proof follows the same steps as in that of Theorem \ref{t2.1} and
the details are omitted.
\end{proof}

\begin{remark}
If $\left\Vert a\right\Vert =1$ and $\left\Vert x_{j}-a\right\Vert \leq 1,$
then (\ref{2.10}) has a simpler form:%
\begin{equation}
\frac{\left\Vert \sum_{j=1}^{n}p_{j}x_{j}\right\Vert }{\sum_{j=1}^{n}p_{j}%
\left\Vert x_{j}\right\Vert }\geq \frac{1}{2}\min_{1\leq j\leq n}\left\{ 
\frac{1-\left\Vert x_{j}-a\right\Vert ^{2}}{\left\Vert x_{j}\right\Vert }%
\right\} \left( \geq 0\right) ,  \label{2.15}
\end{equation}%
while (\ref{2.14}) becomes%
\begin{equation}
\frac{\left\Vert \sum_{j=1}^{n}p_{j}x_{j}\right\Vert }{\sum_{j=1}^{n}p_{j}%
\left\Vert x_{j}\right\Vert }\geq \frac{1}{2}\left( 1-\max_{1\leq j\leq
n}\left\Vert x_{j}-a\right\Vert \right) \left( \geq 0\right) .  \label{2.16}
\end{equation}
\end{remark}

A different approach for bounding the semi-inner product is incorporated in
the following:

\begin{lemma}
\label{l2.2}Let $\left( X,\left\Vert \cdot \right\Vert \right) $ be a normed
space. If $x,a\in X,$ then%
\begin{equation}
\left\langle x,a\right\rangle _{i}\geq \left\Vert a\right\Vert \left(
\left\Vert a\right\Vert -\left\Vert x-a\right\Vert \right) .  \label{b.1}
\end{equation}%
The inequality (\ref{b.1}) is sharp.
\end{lemma}

\begin{proof}
If $a=0,$ then obviously (\ref{b.1}) holds with equality. For $a\neq 0,$
consider%
\begin{equation*}
\tau _{-}\left( x,a\right) :=\lim_{s\rightarrow 0-}\frac{\left\Vert
a+sx\right\Vert -\left\Vert a\right\Vert }{s}.
\end{equation*}%
Observe that%
\begin{align}
\left\langle x,a\right\rangle _{i}& =\lim_{s\rightarrow 0-}\frac{\left\Vert
a+sx\right\Vert ^{2}-\left\Vert a\right\Vert ^{2}}{2s}  \label{b.2} \\
& =\tau _{-}\left( x,a\right) \lim_{s\rightarrow 0-}\left[ \frac{\left\Vert
a+sx\right\Vert +\left\Vert a\right\Vert }{2}\right] =\tau _{-}\left(
x,a\right) \left\Vert a\right\Vert .  \notag
\end{align}%
On the other hand, since the function $R\ni s\longmapsto \left\Vert
a+sx\right\Vert \in \mathbb{R}_{+}$ is convex on $\mathbb{R}$, hence%
\begin{equation}
\tau _{-}\left( x,a\right) \geq \frac{\left\Vert a+\left( -1\right)
x\right\Vert -\left\Vert a\right\Vert }{\left( -1\right) }=\left\Vert
a\right\Vert -\left\Vert x-a\right\Vert .  \label{b.3}
\end{equation}%
Consequently, by (\ref{b.2}) and (\ref{b.3}) we get (\ref{b.1}).

Now, let $x=\varepsilon a,$ $\varepsilon \in \left( 0,1\right) ,$ $a\neq 0.$
Then%
\begin{equation*}
\left\langle x,a\right\rangle _{i}=\varepsilon \left\Vert a\right\Vert
^{2},\quad \left\Vert a\right\Vert -\left\Vert x-a\right\Vert =\left\Vert
a\right\Vert -\left( 1-\varepsilon \right) \left\Vert a\right\Vert
=\varepsilon \left\Vert a\right\Vert ,
\end{equation*}%
which shows that the equality case in (\ref{b.1}) holds true for the nonzero
quantities $\varepsilon \left\Vert a\right\Vert ^{2}$. The proof is complete.
\end{proof}

The following reverse of the generalised triangle inequality may be stated.

\begin{theorem}
\label{t2.3}Let $a,x_{j}\in X\backslash \left\{ 0\right\} $ for $j\in
\left\{ 1,\dots ,n\right\} $ with the property that $\left\Vert a\right\Vert
\geq \left\Vert x_{j}-a\right\Vert $ for $j\in \left\{ 1,\dots ,n\right\} .$
Then for any $p_{j}\geq 0,j\in \left\{ 1,\dots ,n\right\} $ with $%
\sum_{j=1}^{n}p_{j}=1$ we have%
\begin{equation}
\frac{\left\Vert \sum_{j=1}^{n}p_{j}x_{j}\right\Vert }{\sum_{j=1}^{n}p_{j}%
\left\Vert x_{j}\right\Vert }\geq \min_{1\leq j\leq n}\left\{ \frac{%
\left\Vert a\right\Vert -\left\Vert x_{j}-a\right\Vert }{\left\Vert
x_{j}\right\Vert }\right\} \quad \left( \geq 0\right) .  \label{b.4}
\end{equation}%
The inequality (\ref{b.4}) is sharp.
\end{theorem}

\begin{proof}
On making use of Lemma \ref{l2.2}, we have:%
\begin{align*}
\frac{\left\langle x_{j},a\right\rangle _{i}}{\left\Vert x_{j}\right\Vert }&
\geq \left\Vert a\right\Vert \left( \frac{\left\Vert a\right\Vert
-\left\Vert x_{j}-a\right\Vert }{\left\Vert x_{j}\right\Vert }\right) \\
& \geq \left\Vert a\right\Vert \eta ,
\end{align*}%
for each $j\in \left\{ 1,\dots ,n\right\} ,$ where%
\begin{equation*}
\eta :=\min_{1\leq j\leq n}\left\{ \frac{\left\Vert a\right\Vert -\left\Vert
x_{j}-a\right\Vert }{\left\Vert x_{j}\right\Vert }\right\} .
\end{equation*}%
Now utilising the same argument explained in the proof of Theorem \ref{t2.2}%
, we get the desired inequality (\ref{b.4}).

If we choose in (\ref{b.4}) $x_{1}=\cdots =x_{n}=\varepsilon a,$ $%
\varepsilon \in \left( 0,1\right) ,$ $a\neq 0,$ and $p_{1}=...=p_{n}=1$ then
we have equality, and the proof is complete.
\end{proof}

\begin{remark}
The above result may be stated in a simpler way, i.e., if $\rho \in \left(
0,1\right) ,$ $a$ and $x_{j}\in X\backslash \left\{ 0\right\} ,$ $j\in
\left\{ 1,\dots ,n\right\} $ are such that%
\begin{equation}
\left( \left\Vert x_{j}\right\Vert \geq \right) \left\Vert a\right\Vert
-\left\Vert x_{j}-a\right\Vert \geq \rho \left\Vert x_{j}\right\Vert \quad
\left( \geq 0\right)   \label{b.5}
\end{equation}%
for each $j\in \left\{ 1,\dots ,n\right\} ,$ then%
\begin{equation}
\left\Vert \sum_{j=1}^{n}p_{j}x_{j}\right\Vert \geq \rho
\sum_{j=1}^{n}p_{j}\left\Vert x_{j}\right\Vert .  \label{b.6}
\end{equation}
\end{remark}

\section{Other Related Results to the Triangle Inequality}

The following result may be stated:

\begin{theorem}
\label{t3.1}Let $\left( X,\left\Vert \cdot \right\Vert \right) $ be a normed
linear space and $x_{1},\dots ,x_{n}$ nonzero vectors in $X$ and $p_{j}\geq 0
$ with $\sum_{j=1}^{n}p_{j}=1.$ If $\bar{x}_{p}:=\sum_{j=1}^{n}p_{j}x_{j}%
\neq 0$ and there exists a $r>0$ with%
\begin{equation}
\frac{\left\langle x_{j},\bar{x}_{p}\right\rangle _{i}}{\left\Vert
x_{j}\right\Vert \left\Vert \bar{x}_{p}\right\Vert }\geq r\quad \text{for
each \ }j\in \left\{ 1,\dots ,n\right\}   \label{3.1}
\end{equation}%
then%
\begin{equation}
\left\Vert \sum_{j=1}^{n}p_{j}x_{j}\right\Vert \geq
r\sum_{j=1}^{n}p_{j}\left\Vert x_{j}\right\Vert .  \label{3.2}
\end{equation}%
If $p_{j}>0$ for each $j\in \left\{ 1,\dots ,n\right\} ,$ then the equality
holds in (\ref{3.2}) if and only if the equality case hold in (\ref{3.1})
for each $j\in \left\{ 1,\dots ,n\right\} .$
\end{theorem}

\begin{proof}
From (\ref{3.1}) on multiplying with $p_{i}\geq 0$ we have%
\begin{equation*}
\left\langle p_{j}x_{j},\bar{x}_{p}\right\rangle _{i}\geq rp_{j}\left\Vert 
\bar{x}_{p}\right\Vert \left\Vert x_{j}\right\Vert 
\end{equation*}%
for any $j\in \left\{ 1,\dots ,n\right\} .$

Summing over $j$ from $1$ to $n$ and taking into account the superadditivity
property of the interior semi-inner product, we have%
\begin{equation}
\left\langle \sum_{j=1}^{n}p_{j}x_{j},\bar{x}_{p}\right\rangle _{i}\geq
\sum_{j=1}^{n}\left\langle p_{j}x_{j},\bar{x}_{p}\right\rangle _{i}\geq
r\left\Vert \bar{x}_{p}\right\Vert \sum_{j=1}^{n}p_{j}\left\Vert
x_{j}\right\Vert  \label{3.3}
\end{equation}%
and since%
\begin{equation*}
\left\langle \sum_{j=1}^{n}p_{j}x_{j},\bar{x}_{p}\right\rangle
_{i}=\left\Vert \sum_{j=1}^{n}p_{j}x_{j}\right\Vert ^{2}\neq 0
\end{equation*}%
hence by (\ref{3.3}) we get (\ref{3.2}).

The equality case is obvious and the proof is complete.
\end{proof}

For the system of vectors $x_{1},\dots ,x_{k}\in X,$ we denote by $\bar{x}$
its gravity center, i.e.,%
\begin{equation*}
\bar{x}:=\frac{1}{n}\sum_{i=1}^{n}x_{i}.
\end{equation*}

The following corollary is obvious.

\begin{corollary}
\label{c3.1}Let $x_{1},\dots ,x_{n}\in X\backslash \left\{ 0\right\} $ be
such that $\bar{x}\neq 0.$ If there exists a $r>0$ such that%
\begin{equation}
\frac{\left\langle x_{j},\bar{x}\right\rangle _{i}}{\left\Vert
x_{j}\right\Vert \left\Vert \bar{x}\right\Vert }\geq r\quad \text{for each \ 
}j\in \left\{ 1,\dots ,n\right\} ,  \label{3.4}
\end{equation}%
then the following reverse of the generalised triangle inequality holds:%
\begin{equation}
\left\Vert \sum_{j=1}^{n}x_{j}\right\Vert \geq r\sum_{j=1}^{n}\left\Vert
x_{j}\right\Vert .  \label{3.5}
\end{equation}%
The equality holds in (\ref{3.5}) if and only if the case of equality holds
in (\ref{3.4}) for each $j\in \left\{ 1,\dots ,n\right\} .$
\end{corollary}

The following refinements of the generalised triangle inequality may be
stated as well:

\begin{theorem}
\label{t3.2}Let $x_{i},\bar{x}_{p},p_{i},$ $i\in \left\{ 1,\dots ,n\right\} $
be as in Theorem \ref{t3.1}. If there exists a constant $R$ with $1>R>0$ and
such that%
\begin{equation}
R\geq \frac{\left\langle x_{j},\bar{x}_{p}\right\rangle _{s}}{\left\Vert
x_{j}\right\Vert \left\Vert \bar{x}_{p}\right\Vert }\quad \text{for each \ }%
j\in \left\{ 1,\dots ,n\right\} ,  \label{3.6}
\end{equation}%
then%
\begin{equation}
R\sum_{j=1}^{n}p_{j}\left\Vert x_{j}\right\Vert \geq \left\Vert
\sum_{j=1}^{n}p_{j}x_{j}\right\Vert .  \label{3.7}
\end{equation}%
If $p_{j}>0$ for each $j\in \left\{ 1,\dots ,n\right\} ,$ then the equality
holds in (\ref{3.7}) if and only if the equality case holds in (\ref{3.6})
for each $j\in \left\{ 1,\dots ,n\right\} .$
\end{theorem}

The proof is similar to the one in Theorem \ref{t3.1} on taking into account
that the superior semi-inner product is a subadditative functional in the
first variable.

\begin{corollary}
\label{c3.2}Let $x_{j},$ $j\in \left\{ 1,\dots ,n\right\} $ be as in
Corollary \ref{c3.1}. If there exists an $R$ with $1>R>0$ and%
\begin{equation}
R\geq \frac{\left\langle x_{j},\bar{x}\right\rangle _{s}}{\left\Vert
x_{j}\right\Vert \left\Vert \bar{x}\right\Vert }\quad \text{for each \ }j\in
\left\{ 1,\dots ,n\right\} ,  \label{3.8}
\end{equation}%
then the following refinement of the generalised triangle inequality holds:%
\begin{equation}
R\sum_{j=1}^{n}\left\Vert x_{j}\right\Vert \geq \left\Vert
\sum_{j=1}^{n}x_{j}\right\Vert .  \label{3.9}
\end{equation}%
The equality hold in (\ref{3.9}) if and only if the case of equality holds
in (\ref{3.8}) for each $j\in \left\{ 1,\dots ,n\right\} .$
\end{corollary}

\end{document}